\documentclass[11pt,reqno]{amsart}
\usepackage{amsmath}
\usepackage{amssymb}
\usepackage[bookmarks]{hyperref}
\newtheorem{theorem}{Theorem}[section]

\newtheorem{corollary}[theorem]{Corollary}

\newtheorem{remark}[theorem]{Remark}

\numberwithin{equation}{section}

\allowdisplaybreaks[1]
\begin{document}

\title[Finite-rank products of Toeplitz operators]{A refined Luecking's theorem and finite-rank products of Toeplitz operators}
\author{Trieu Le}
\address{Trieu Le, Department of Mathematics, University of Toronto, Toronto, Ontario, Canada M5S 2E4}
\email{trieu.le@utoronto.ca}
\begin{abstract} For any function $f$ in $L^{\infty}(\mathbb{D})$, let $T_f$ denote the corresponding Toeplitz operator the Bergman space $A^2(\mathbb{D})$. A recent result of D. Luecking shows that if $T_f$ has finite rank then $f$ must be the zero function. Using a refined version of this result, we show that if all except possibly one of the functions $f_1,\ldots, f_{m}$ are radial and $T_{f_1}\cdots T_{f_m}$ has finite rank, then one of these functions must be zero.
\end{abstract}
\maketitle

\section{Introduction}
As usual, let $\mathbb{D}$ denote the unit disk and $\mathbb{T}$ denote the unit circle in the complex plane $\mathbb{C}$. Let $\mathrm{d}A$ denote the Lebesgue measure on $\mathbb{D}$ which is normalized such that the unit disk has total mass $1$. We have $\mathrm{d}A(z)=\frac{1}{\pi}\mathrm{d}x\mathrm{d}y$, where $z=x+\mathrm{i}y$ for $x,y$ real. We write $L^2$ for $L^2(\mathbb{D},\mathrm{d}A)$. The Bergman space $A^2$ is the subspace of $L^2$ that consists of holomorphic functions. It is well-known that $A^2$ is a closed subspace of $L^2$. The standard orthonormal basis for $A^2$ is $\{e_{m}: m=0,1,\ldots\}$, where $e_{m}(z)=\sqrt{m+1}\ z^m$ for any non-negative integer $m$. Let $P$ denote the orthogonal projection from $L^2$ onto $A^2$. For any function  $f\in L^2$, the Toeplitz operator with symbol $f$ is denoted by $T_f$, which is densely defined on $A^2$ by $T_f\varphi = P(f\varphi)$ for $\varphi\in H^{\infty}$ - the space of all bounded holomorphic functions on $\mathbb{D}$.
The operator $T_f$ is in fact an integral operator by the formula
\begin{equation*}
(T_f\varphi)(z) = \int\limits_{\mathbb{D}}\dfrac{f(w)\varphi(w)}{(1-\bar{w}z)^{2}}\mathrm{d}A(w),\text{ for } z\in\mathbb{D}, \varphi\in H^{\infty}.
\end{equation*}
If $f$ is a bounded function then $T_{f}$ is a bounded operator on $A^{2}$ with $\|T_f\|\leq\|f\|_{\infty}$ and $(T_{f})^{*}=T_{\bar{f}}$. However, unbounded symbol can also give rise to bounded Toeplitz operators. In fact, since $T_f$ is an integral operator with kernel $\dfrac{f(w)}{(1-\bar{w}z)^2}$ for $z,w\in\mathbb{D}$, we see that if $f\in L^2$ supported in a compact subset of $\mathbb{D}$ then $T_f$ is a compact operator on $A^2$. 

A function $f$ on $\mathbb{D}$ is called a radial function if we have $f(z)=f(|z|)$ for almost all $z\in\mathbb{D}$. If $f\in L^2$ is radial then using polar coordinate we see that
\begin{align*}
\langle T_{f}e_{m}, e_{k}\rangle & = \sqrt{(m+1)(k+1)}\int_{\mathbb{D}}f(z)z^{m}\bar{z}^{k}\mathrm{d}A(z)\\
& = \begin{cases}
0 & \text{if } m\neq k\\
(m+1)\int_{0}^{1}2f(t)t^{2m+1}\mathrm{d}t & \text{if } m=k
\end{cases}\\
& = \begin{cases}
0 & \text{if } m\neq k\\
(m+1)\int_{0}^{1}f(r^{1/2})r^{m}\mathrm{d}r & \text{if } m=k.
\end{cases}
\end{align*}
This shows that the operator $T_f$ is diagonal with respect to the standard orthonormal basis. The eigenvalues of $T_f$ are given by
\begin{equation}\label{eqn-500}
\omega(f,m) = \langle T_fe_{m}, e_{m}\rangle = (m+1)\int_{0}^{1}f(r^{1/2})r^{m}\mathrm{d}r,\quad m=0,1,\ldots.
\end{equation}

It follows from Stone-Weierstrass's Theorem that if $f\in L^2$ such that $T_f$ is the zero operator then $f$ must vanish almost everywhere in $\mathbb{D}$. On the other hand, the problem of determining whether there exists a nontrivial finite rank Toeplitz operator on $A^2$ was open for quite a long time. Recently D. Luecking has found an elegant proof that gives the negative answer to this problem.

There is an extensive literature on Toeplitz operators on the Hardy space $H^2$ of the unit circle. We refer the reader to \cite{Martinez-Avendano2007} for definitions of $H^2$ and their Toeplitz operators. In the context of Toeplitz operators on $H^2$, it was showed by A. Brown and P. Halmos \cite{Brown1963} back in the sixties that if $f$ and $g$ are bounded functions on the unit circle then $T_{g}T_{f}$ is another Toeplitz operator if and only if either $f$ or $\bar{g}$ is holomorphic. From this it is readily deduced that if $f,g\in L^{\infty}(\mathbb{T})$ such that $T_{g}T_{f}=0$ then one of the symbols must be the zero function. In contrast with this, for Toeplitz operators on the Bergman space, it has not been known if it is true that for $f,g\in L^{\infty}(\mathbb{D})$, $T_{g}T_{f}=0$ implies $g$ or $f$ is the zero function. Affirmative answers have been obtained by researchers only in special cases. In \cite{Ahern2001}, P. Ahern and {\v{Z}}. {\v{C}}u{\v{c}}kovi{\'c} answered this problem affirmatively with the assumption that both $f$ and $g$ are bounded harmonic functions on $\mathbb{D}$. Later in \cite{Cuckovic2003}, {\v{C}}u{\v{c}}kovi{\'c} was able to show that if $f,g$ are bounded such that $f$ is harmonic and $g(r\mathrm{e}^{\mathrm{i}\theta})=\sum_{m=-\infty}^{N}g_{k}(r)\mathrm{e}^{\mathrm{i}m\theta}$ for $z=r\mathrm{e}^{\mathrm{i}\theta}\in\mathbb{D}$, then $T_{g}T_{f}=0$ implies either $f=0$ or $g=0$. The case one of the symbols is a bounded radial function has also been understood. See \cite{Ahern2004} and \cite{Le-5} for more details. In fact, in \cite{Le-5}, the author was able to show that if all except possibly one of the functions $f_1,\ldots, f_M$ are bounded radial functions and $T_{f_1}\cdots T_{f_M}=0$ then one of these functions must be zero.

A more general problem than the above zero product problem is the finite rank product problem. Recall that the above mentioned theorem of Luecking shows that if $f\in L^{2}$ such that $T_f$ has finite rank then $f$ is the zero function. What happens if $T_gT_f$ has finite rank, where $f$ and $g$ are bounded measurable functions on the unit disk? The answer to this general question seems to be still far from completed but the following important case has been understood: If $f$ and $g$ are bounded harmonic functions then one of them must be the zero function (K. Guo, S. Sun and D. Zheng \cite{Guo2007}). The purpose of this paper is to report the same answer in some other special cases.

In the first part of this paper, we use Luecking's Theorem to show that if $f,g$ are functions in $L^2$ where $f$ satisfies a certain condition and $T_gT_f$ (which is densely defined on $A^2$) has finite rank, then either $f=0$ or $g=0$. In the second part of the paper, we prove a ``refined'' version of Luecking's Theorem and use it to show that if $f_{1},\ldots, f_{m_1}$ and $g_{1},\ldots,g_{m_2}$ are radial functions in $L^{\infty}$ and $f$ is a function in $L^2$ such that $T_{g_1}\cdots T_{g_{m_2}}T_fT_{f_1}\cdots T_{f_{m_1}}$ (which is densely defined on $A^2$) has finite rank, then one of the above functions must be zero.

\section{Finite rank products of two Toeplitz operators}\label{section-2}
We begin this section with a detailed discussion of the decomposition $L^2=\bigoplus_{m\in\mathbb{Z}}\mathcal{R}\mathrm{e}^{\mathrm{i}m\theta}$, where
\begin{equation*}
\mathcal{R}=\{u:[0,1)\rightarrow\mathbb{C}\text{ such that }\int_{0}^{1}|u(r)|^2r\mathrm{d}r<\infty\}.
\end{equation*}
This decomposition has been used by {\v{C}}u{\v{c}}kovi{\'c} and Rao in their studies of Toeplitz operators (see Section 2 in \cite{Cuckovic1998}). Let $f\in L^{2}(\mathbb{D})$. Then for almost all $r\in [0,1)$, the function $\zeta\mapsto f(r\zeta)$ for $\zeta\in\mathbb{T}$ is in $L^2(\mathbb{T},\frac{1}{2\pi}\mathrm{d}\theta)$. Since $\{\zeta^{m}: m\in\mathbb{Z}\}$ is an orthonormal basis for $L^2(\mathbb{T},\frac{1}{2\pi}\mathrm{d}\theta)$, we have
\begin{equation*}
f(r\zeta) = \sum\limits_{m=-\infty}^{\infty}\Big(\dfrac{1}{2\pi}\int\limits_{0}^{2\pi}f(r\mathrm{e}^{\mathrm{i}\theta})\mathrm{e}^{-\mathrm{i}m\theta}\mathrm{d}\theta\Big)\zeta^{m},
\end{equation*}
where the sum takes place in $L^2(\mathbb{T})$. For $m\in\mathbb{Z}$, define
\begin{equation*}
f_{m}(r) = \dfrac{1}{2\pi}\int\limits_{0}^{2\pi}f(r\mathrm{e}^{-\mathrm{i}m\theta})\mathrm{d}\theta,\quad 0\leq r <1.
\end{equation*}
Then the above representation becomes (with $\zeta=\mathrm{e}^{\mathrm{i}\theta}$),
\begin{equation}\label{eqn-1}
f(r\mathrm{e}^{\mathrm{i}\theta}) = \sum\limits_{m=-\infty}^{\infty} f_{m}(r)\mathrm{e}^{\mathrm{i}m\theta}.
\end{equation}
This representation holds for almost all $r\in [0,1)$ and for such $r$, the sum on the right hand side takes place in $L^2(\mathbb{T})$. Now we have
\begin{align*}
\|f\|^2_{L^2(\mathbb{D})} & = \int\limits_{0}^{1}\Big(\dfrac{1}{2\pi}\int\limits_{0}^{2\pi}|f(r\mathrm{e}^{\mathrm{i}\theta})|^2\mathrm{d}\theta\Big)r\mathrm{d}r\\
& = \int\limits_{0}^{1}\Big(\sum\limits_{m=-\infty}^{\infty}|f_{m}(r)|^2\Big)r\mathrm{d}r\\
& = \sum\limits_{m=-\infty}^{\infty}\int\limits_{0}^{1}|f_{m}(r)|^2r\mathrm{d}r.
\end{align*}
This shows that $f_m\in\mathcal{R}$ for all $m\in\mathbb{Z}$ and the right hand side of \eqref{eqn-1} converges in $L^2(\mathbb{D})$. Therefore the representation \eqref{eqn-1} in fact takes place in $L^2(\mathbb{D})$.

The following theorem is our first result in the paper.

\begin{theorem}\label{theorem-1} Suppose $f\in L^{2}$ with $f(r\mathrm{e}^{\mathrm{i}\theta}) =\!\!\! \sum\limits_{m=-\infty}^{M} f_{m}(r)\mathrm{e}^{\mathrm{i}m\theta}$ for $z=r\mathrm{e}^{\mathrm{i}\theta}$, where $M$ is an integer. Assume that $\int_{0}^{1}f_{M}(r)r^{k}\mathrm{d}r\neq 0$ for all $k\geq N$, where $N$ is a positive integer. If $g\in L^{2}$ such that $T_gT_f$ (which is densely defined on $A^2$) has finite rank then $g$ is the zero function.
\end{theorem}

\begin{proof} Recall that $A^2(\mathbb{D})$ has the orthonormal basis $\{e_{m}: m=0,1,\ldots\}$, where $e_{m}(z)=\sqrt{m+1}\ z^m$ for any non-negative integer $m$. For any non-negative integers $k,l$ we have
\begin{align*}
\langle T_f e_{k}, e_{l}\rangle & = \sqrt{(k+1)(l+1)}\int\limits_{\mathbb{D}}f(z)z^{k}\bar{z}^l\mathrm{d}A(z)\\
& = \sqrt{(k+1)(l+1)}\int\limits_{0}^{1}\Big(\dfrac{1}{2\pi}f(r\mathrm{e}^{\mathrm{i}\theta})\mathrm{e}^{\mathrm{i}(k-l)\theta}\mathrm{d}\theta\Big)r^{k+l+1}\mathrm{d}r\\
& = \sqrt{(k+1)(l+1)}\int\limits_{0}^{1}f_{l-k}(r)r^{k+l+1}\mathrm{d}r.
\end{align*}
By assumption about $f$, $\langle T_f e_{k}, e_{l}\rangle = 0$ whenever $l-k>M$. Thus for $k\in\mathbb{N}$ such that $k+M\geq 0$, we have
\begin{align*}
T_fe_{k} & = \sum\limits_{l=0}^{\infty}\langle T_f e_{k}, e_{l}\rangle e_{l}\\
& = \sqrt{k+1}\sum\limits_{l=0}^{k+M}\Big(\sqrt{l+1}\int\limits_{0}^{1}f_{l-k}(r)^{k+l+1}\mathrm{d}r\Big)e_{l}\\
& = \sqrt{(k+1)(M+k+1)}\Big(\int\limits_{0}^{1}f_{M}(r)r^{2k+M+1}\mathrm{d}r\Big)e_{k+M}\\
&\quad + \sqrt{k+1}\sum\limits_{l=0}^{k+M-1}\Big(\sqrt{l+1}\int\limits_{0}^{1}f_{l-k}(r)r^{k+l+1}\mathrm{d}r\Big)e_{l}
\end{align*}
This shows that when $k+M\geq 1$ and $2k+M+1\geq N$, $e_{k+M}$ can be written as a linear combination of $\{T_{f}e_{k}\}\cup\{e_{0},\ldots, e_{k+M-1}\}$.

Now suppose $T_gT_f$ has finite rank and let $\{\varphi_1,\ldots,\varphi_K\}$ is a set that spans $T_gT_f(\mathcal{P})$ where $\mathcal{P}$ is the space of all polynomials in the variable $z$. Then for any non-negative integer $k$ with $k+M\geq 1$ and $2k+M+1\geq N$ we see that $T_ge_{k+M}$ is a linear combination of $\{\varphi_{1}, \ldots, \varphi_{K}\}\cup\{T_g(e_{0}),\ldots, T_{g}(e_{k+M-1})\}$. From this, it follows by induction that $T_g$ is a finite rank operator. By Luecking's Theorem \cite{Luecking2008} or a refined version of it (see Theorem \ref{theorem-2} in Section \ref{section-3}), we see that $g$ is the zero function.
\end{proof}

\begin{remark} If $f(z)=\bar{h}(z)+p(z,\bar{z})$ where $h\in A^2$ and $p$ a polynomial in two variables then $f$ can be written in the form in the hypothesis of Theorem \ref{theorem-1}. Therefore, Theorem \ref{theorem-1} shows that if $T_gT_f$ is of finite rank for some $g\in L^2$ then either $f$ or $g$ is the zero function.
\end{remark}

\section{A refined Luecking's Theorem and finite rank products of Toeplitz operators}\label{section-3}

We begin this section by a refined version of Luecking's Theorem whose proof is greatly influenced by Luecking's argument. For the rest of the paper, let $\mathcal{P}$ denote the space of all polynomials in the variable $z$.

\begin{theorem}\label{theorem-2} Let $\mathcal{S}\subset\mathbb{N}$ ($\mathbb{N}$ denotes the set of all non-negative integers) so that $\sum_{s\in\mathcal{S}}\frac{1}{s+1}<\infty$. Let $\mathcal{N}$ be the subspace of $\mathcal{P}$ spanned by the monomials $\{z^{m}:m\in\mathbb{N}\backslash\mathcal{S}\}$ and let $\mathcal{N}^{*}=\{\bar{g}: g\in\mathcal{N}\}$. Let $\nu$ be a complex regular Borel measure on $\mathbb{C}$ with compact support. Let $T_{\nu}$ be the operator from $\mathcal{N}$ to the space of linear functionals on $\mathcal{N}^{*}$ by $T_{\nu}f(\bar{g})=\int_{\mathbb{C}}f\bar{g}\mathrm{d}\nu$ for all $f,g\in\mathcal{N}$. Then $T_{\nu}$ has finite rank if and only if the support of $\nu$ is finite.
\end{theorem}

\begin{proof}
For any $z\in\mathbb{C}$, let $\delta_{z}$ denote the point mass measure concentrated at $z$. Since $T_{\nu-\nu(\{0\})\delta_{0}} = T_{\nu}-\nu(\{0\})T_{\delta_{0}}$, we see that $T_{\nu}$ has finite rank if and only if $T_{\nu-\nu(\{0\})\delta_{0}}$ has finite rank. So without loss of generality, we may assume that $\nu(\{0\})=0$.

If the support of $\nu$ is contained in a finite set $\{z_1,\ldots, z_{N-1}\}$ for some $N\geq 2$, then $T_{\nu}=\sum_{j=1}^{N-1}\nu(\{z_j\})T_{\delta_{z_j}}$. Hence $T_{\nu}$ has rank less than $N$.

Conversely, suppose $T_{\nu}$ has rank less than $N$. Following Luecking's argument in \cite[p. 3]{Luecking2008}, we see that for any $f_1,\ldots,f_N$ and $g_1,\ldots,g_N$ in $\mathcal{N}$,
\begin{equation}\label{eqn-400}
\int_{\mathbb{C}^n}\prod_{l=1}^{N}f_{l}(z_l)\det(\bar{g}_i(z_j))\mathrm{d}\nu^{N}(Z)=0,
\end{equation}
where $Z=(z_1,\ldots,z_N)\in\mathbb{C}^N$ and $\nu^{N}$ is the product of $N$ copies of $\mu$ on $\mathbb{C}^{N}$.

Let $m_1,\ldots,m_N$ and $k_1,\ldots,k_N$ be non-negative integers. Let
\begin{align*}
\mathcal{Z} & = \{s\in\mathbb{N}: s+m_j\notin\mathcal{S} \text{ and } s+k_j\notin\mathcal{S} \text{ for all } 1\leq j\leq N\}\\
& = \mathbb{N}\backslash\Big((\cup_{j=1}^{N}(\mathcal{S}-m_j))\cup(\cup_{j=1}^{N}\mathcal({S}-k_j))\Big).
\end{align*}
Since $\sum\limits_{s\in\mathcal{S}}\frac{1}{s+1}<\infty$ we have $\sum\limits_{s\in\mathbb{N}\backslash\mathcal{Z}}\frac{1}{s+1}<\infty$. This shows that
\begin{equation}\label{eqn-401}
\sum\limits_{s\in\mathcal{Z}}\frac{1}{s+1} = \infty.
\end{equation}
Now for any $s\in\mathcal{Z}$, the monomials $f_{j}(z)=z^{m_j+s}$ and $g_{j}(z)=z^{k_j+s}$ for $j=1,\ldots,N$ are not in $\mathcal{N}$. So we may use \eqref{eqn-400} to get
\begin{align}
0 & = \int_{\mathbb{C}^n}\prod_{l=1}^{N}z_{l}^{m_l+s}\det(\bar{z}_j^{k_i+s})\mathrm{d}\nu^{N}(Z)\notag\\
& = \int_{\mathbb{C}^n}\prod_{l=1}^{N}z_{l}^{m_l}\det(\bar{z}_j^{k_i})|z_{1}\ldots z_{N}|^{2s}\mathrm{d}\nu^{N}(Z)\notag\\
& = \int_{\mathbb{C}^n\backslash W}\prod_{l=1}^{N}z_{l}^{m_l}\det(\bar{z}_j^{k_i})|z_{1}\ldots z_{N}|^{2s}\mathrm{d}\nu^{N}(Z),\label{eqn-402}
\end{align}
where $W=\{Z=(z_1,\ldots,z_N)\in\mathbb{C}^N: z_1\cdots z_N=0\}$. The last identity follows from the fact that $\nu^N(W)=0$.

Let $\mathbb{K}$ denote the open right half plane consisting of all $w$ with $\Re(w)>0$ and let $\bar{\mathbb{K}}$ denote the closure of $\mathbb{K}$ in $\mathbb{C}$. For any $w\in\mathbb{K}$, define
\begin{equation*}
F(w) = \int_{\mathbb{C}^n\backslash W}\prod_{l=1}^{N}z_{l}^{m_l}\det(\bar{z}_j^{k_i})|z_{1}\ldots z_{N}|^{2w}\mathrm{d}\nu^{N}(Z).
\end{equation*}
Here, for a positive number $t$ and a complex number $w$, $t^{w}=\exp(w\log t)$ where $\log$ is the principal branch of the logarithmic function.

Suppose the measure $\nu$ is supported in the disk $D(0,R)$ of radius $R>0$ centered at the origin in the complex plane. Then $\nu^{N}$ is supported in the polydisk $D_{N}(0,R)$ of the same radius centered at the origin in $\mathbb{C}^N$. Then for any $w\in\mathbb{K}$ and any $Z=(z_1,\ldots,z_N)$ in the above polydisk, we have
\begin{equation*}
\big||z_1\cdots z_N|^{2w}\big| = |z_1\cdots z_N|^{2\Re(w)}\leq R^{2N\Re(w)}.
\end{equation*}
Therefore,
\begin{equation*}
|F(w)| = \Big|\int_{D_{N}(0,R)\backslash W}\prod_{l=1}^{N}z_{l}^{m_l}\det(\bar{z}_j^{k_i})|z_{1}\cdots z_{N}|^{2w}\mathrm{d}\nu^{N}(Z)\Big|\leq CR^{2N\Re(w)},
\end{equation*}
where $C$ is a constant independent of $w$. It follows that $F$ is not only defined on $\bar{\mathbb{K}}$ but also continuous on $\bar{\mathbb{K}}$. Now an application of Morera's Theorem shows that $F$ is analytic on $\mathbb{K}$. Let $G(w)=F(w)R^{-2Nw}$ for $w\in\mathbb{K}$, then $G$ is continuous, bounded on $\bar{\mathbb{K}}$ and analytic on $\mathbb{K}$. Now define
\begin{equation*} H(\zeta) = G\Big(\dfrac{1+\zeta}{1-\zeta}\Big)\quad\quad (|\zeta|<1).\end{equation*}
Then $H$ is a bounded analytic function on the unit disk. For any $s\in\mathcal{Z}$, equation \eqref{eqn-402} and the definitions of $F,G$ show that $G(s)=F(s)=0$, which implies $H(\frac{s-1}{s+1})=0$. Now
\begin{equation*}
\sum\limits_{\substack{s\in\mathcal{Z}\\ s\geq 1}}(1-|\dfrac{s-1}{s+1}|) = \sum\limits_{\substack{s\in\mathcal{Z}\\ s\geq 1}}\dfrac{2}{s+1} = \infty\quad\text{ by \eqref{eqn-401}}.
\end{equation*}
Corollary to Theorem 15.23 in \cite{Rudin1987} shows that $H$ is identically zero on the unit disk. Hence $G$ and $F$ are identically zero in $\bar{\mathbb{K}}$. In particular, $F(0)=0$, which shows that
\begin{equation*}
\int_{\mathbb{C}^n}\prod_{l=1}^{N}z_{l}^{m_l}\det(\bar{z}_j^{k_i})\mathrm{d}\nu^{N}(Z) = 0.
\end{equation*}
Since $m_1,\ldots,m_N$ and $k_1,\ldots,k_N$ were arbitrary non-negative integers, we conclude that \eqref{eqn-400} holds for all $f_1,\ldots,f_N$ and $g_1,\ldots,g_N$ in $\mathcal{P}$. Following Luecking's argument again \cite[Section 4 and 5]{Luecking2008}, we conclude that the support of $\nu$ is finite.
\end{proof}

Now let $\mathcal{S}$ and $\mathcal{N}$ be as in the hypothesis of Theorem \ref{theorem-2}. Let $\mathcal{M}$ denote the subspace of $\mathcal{P}$ spanned by $\{z^{m}: m\in\mathcal{S}\}$. Let $\bar{\mathcal{M}}$ (respectively, $\bar{\mathcal{N}}$) denote the closure of $\mathcal{M}$ (respectively, $\mathcal{N}$) in $A^2$.

\begin{corollary}\label{cor-1} Suppose $f\in L^2$ so the operator $T_{f}$ is densely defined on $A^2$. If $T_{f}(\mathcal{N})\subset \mathrm{Span}(\bar{\mathcal{M}}\cup\{\varphi_{1},\ldots,\varphi_{N}\})$, where $\varphi_{1},\ldots,\varphi_{N}\in A^2$, then $f$ is the zero function.
\end{corollary}

\begin{proof}
Let $P_{\bar{\mathcal{M}}}$ (respectively, $P_{\bar{\mathcal{N}}}$) denote the orthogonal projection from $A^2$ onto $\bar{\mathcal{M}}$ (respectively, $\bar{\mathcal{N}}$). Then we have $P_{\bar{\mathcal{N}}}=1-P_{\bar{\mathcal{M}}}$ and hence $P_{\bar{\mathcal{M}}}P_{\bar{\mathcal{N}}}=P_{\bar{\mathcal{N}}}P_{\bar{\mathcal{M}}}=0$. By replacing $\varphi_j$ by $\varphi_j-P_{\bar{\mathcal{M}}}\varphi_j$ if necessary, we may assume that $\varphi_j\perp\mathcal{M}$ for $1\leq j\leq N$. By using the Gram-Schmidt process if necessary, we may assume that the vectors $\varphi_{1},\ldots,\varphi_{N}$ are orthonormal (we may have fewer vectors after using Gram-Schmidt process but let us still denote by $N$ the total number of the vectors).  

For any $p$ in $\mathcal{N}$ we have $
T_{f}p = P_{\bar{\mathcal{M}}}T_{f}p + \sum\limits_{j=1}^{N}\langle T_{f}p,\varphi_{j}\rangle \varphi_{j}$. This shows that $P_{\bar{\mathcal{N}}}(T_{f}p) = \sum\limits_{j=1}^{N}\langle T_{f}p,\varphi_{j}\rangle P_{\bar{\mathcal{N}}}\varphi_{j} = \sum\limits_{j=1}^{N}\langle T_{f}p,\varphi_{j}\rangle\varphi_{j}=\sum\limits_{j=1}^{N}\langle fp,\varphi_{j}\rangle\varphi_{j}$. Then for any $q$ in $\mathcal{N}$, we have
\begin{align*}
\int_{\mathbb{D}}fp\bar{q}\ \mathrm{d}A & = \langle T_{f}p,q\rangle
 = \langle P_{\bar{\mathcal{N}}}(T_{f}p), q\rangle
 = \sum\limits_{j=1}^{N}\langle fp,\varphi_{j}\rangle\langle\varphi_{j},q\rangle.
\end{align*}
Let $\mathrm{d}\nu = f\mathrm{d}A$. Then the map $T_{\nu}$ from $\mathcal{N}$ to the space of linear functionals on $\mathcal{N}^{*}$ defined by $T_{\nu}p(\bar{q})=\int_{\mathbb{D}}p\bar{q}\mathrm{d}\nu=\int_{\mathbb{D}}fp\bar{q}\mathrm{d}A$ for $p,q\in\mathcal{N}$ is of finite rank. Now Theorem \ref{theorem-2} shows that the support of $\nu$ is finite, which implies that $f(z)=0$ for almost all $z\in\mathbb{D}$.
\end{proof}

\begin{theorem}\label{theorem-3} Suppose $f_{1},\ldots, f_{m_1}$ and $g_1,\ldots,g_{m_2}$ are radial functions in $L^{\infty}$ none of which is the zero function. Suppose $f$ is a function in $L^{2}$ such that $T_{g_1}\cdots T_{g_{m_2}}T_{f}T_{f_1}\cdots T_{f_{m_1}}$ (which is densely defined on $A^2$) is of finite rank, then $f$ must be the zero function.
\end{theorem}

\begin{proof} For any $h\in\{f_1,\ldots, f_{m_1}, g_1,\ldots, g_{m_2}\}$, the operator $T_h$ is diagonal with eigenvalues $\omega(h,m)$ given by \eqref{eqn-500} for $m=0,1,\ldots$. Let $Z(h)=\{m\in\mathbb{N}:\omega(h,m)=0\}$. Since $h$ is not the zero function, {M}\"untz-{S}z\'asz's Theorem (see \cite[Theorem 15.26]{Rudin1987}) shows that $\sum\limits_{m\in Z(h)}\frac{1}{m+1}<\infty$.

Let $\mathcal{S}=Z(f_1)\cup\cdots\cup Z(f_{m_1})\cup Z(g_1)\cup\cdots\cup Z(g_{m_2})$. Then we have $\sum\limits_{m\in\mathcal{S}}\frac{1}{s+1}<\infty$. Let $\mathcal{N}$ (respectively, $\mathcal{M}$) is the subspace of $\mathcal{P}$ spanned by $\{e_{m}: m\in\mathbb{N}\backslash\mathcal{S}\}$ (respectively, $\{e_{m}: m\in\mathcal{S}\}$). Recall that $\mathcal{P}$ denotes the space of all analytic polynomials in the variable $z$.

Put $S_1=T_{f_1}\cdots T_{f_{m_1}}$ and $S_2=T_{g_1}\cdots T_{g_{m_2}}$. For $\varphi\in A^2$ we have
\begin{align*}
S_2\varphi & = T_{g_1}\cdots T_{g_{m_2}}\big(\sum_{j=1}^{\infty}\langle\varphi, e_{j}\rangle e_{j}\big) = \sum_{j=1}^{\infty}\omega(g_1,j)\cdots\omega(g_{m_2},j)\langle\varphi, e_{j}\rangle e_{j}.
\end{align*}
Hence if $S_2\varphi = 0$, then $\omega(g_1,j)\cdots\omega(g_{m_2},j)\langle\varphi, e_{j}\rangle=0$ for all $j\in\mathbb{N}$. It then implies that $\langle\varphi,e_{j}\rangle=0$ whenever $j\in\mathbb{N}\backslash\mathcal{S}$. Thus $\ker(S_2)\subset\bar{\mathcal{M}}$.

On the other hand, if $j\in\mathbb{N}\backslash\mathcal{S}$ then $\omega(f_1,j)\cdots\omega(f_{m_1},j)\neq 0$, and hence,
\begin{equation*}
e_{j} = \dfrac{1}{\omega(f_1,j)\cdots\omega(f_{m_1},j)}T_{f_{1}}\cdots T_{f_{m_1}}e_{j} = \dfrac{1}{\omega(f_1,j)\cdots\omega(f_{m_1},j)}S_1e_{j}.
\end{equation*}
This shows that $\mathcal{N}\subset S_1(\mathcal{N})\subset S_1(\mathcal{P})$. Hence the domain of the operator $S_2T_{f}S_1$ contains $\mathcal{P}$, which is dense in $A^2$.

Now suppose that $S_2T_{f}S_1(\mathcal{P})$ is of finite dimensions, spanned by the set $\{u_1,\ldots, u_N\}$. Let $v_j\in A^2$ such that $S_2v_j = u_j$ for $j=1,\ldots, N$. It then follows that $T_{f}S_1(\mathcal{P})$ is contained in $\mathrm{Span}(\ker(S_2)\cup\{v_1,\ldots, v_N\})$, which is a subspace of $\mathrm{Span}(\bar{\mathcal{M}}\cup\{v_1,\ldots,v_N\})$. But as we have seen above, $\mathcal{N}$ is a subspace of $S_1(\mathcal{P})$. So we conclude that $T_{f}(\mathcal{N})\subset\mathrm{Span}(\bar{\mathcal{M}}\cup\{v_1,\ldots, v_N\})$. Corollary \ref{cor-1} then implies that $f$ is the zero function.
\end{proof}

\begin{remark} Suppose $\mathcal{S}\subset\mathbb{N}$ such that $\sum_{s\in\mathcal{S}}\frac{1}{s+1}<\infty$. Let $\mathcal{N}$ (respectively, $\mathcal{M}$) is the subspace of $\mathcal{P}$ spanned by $\{e_{m}: m\in\mathbb{N}\backslash\mathcal{S}\}$ (respectively, $\{e_{m}: m\in\mathcal{S}\}$). From the proof of Theorem \ref{theorem-3}, we see that if $S_1, S_2$ are bounded operators on $A^2$ such that $\mathcal{N}\subset S_1(\mathcal{P})$, $\ker(S_2)\subset\bar{\mathcal{M}}$ and $S_2T_{f}S_1$ has finite rank then $f$ must be zero. This shows that the conclusion of Theorem \ref{theorem-3} remains valid if $f_{j}(r\mathrm{e}^{\mathrm{i}\theta})=\tilde{f_j}(r)\mathrm{e}^{\mathrm{i}s_j\theta}$ and $g_{k}(r\mathrm{e}^{\mathrm{i}\theta})=\tilde{g_k}(r)\mathrm{e}^{\mathrm{i}t_k\theta}$ for bounded functions $\tilde{f_j},\tilde{g_k}$ on $[0,1)$ and integers $s_j, t_k$, for $1\leq j\leq m_1, 1\leq k\leq m_2$.
\end{remark}

\end{document}